\documentclass[12pt]{article}
\usepackage{amsmath}
\usepackage{amssymb}
\usepackage{vatola}

\def\q{\quad}
\def\qq{\qquad}
\def\mod#1{\ (\text{\rm mod}\ #1)}

\def\t{\hbox}
\def\qtq#1{\q\t{#1}\q}

\def\f{\frac}
\def\e{\equiv}
\def\a{\alpha}
\def\b{\binom}

\def\o{\omega}

 \def\abc#1{$(\t{\rm#1})$}
 \def\an{\{a_n\}}
  \def\ncdots{n=0,1,2,\ldots}
 \def\lra{\Longleftrightarrow}
\def\n{\sum_{k=0}^n\binom nk(-1)^k}

\def\nm{\sum_{k=0}^n\b{n-m}k(-1)^{n-k}}
\def\({\left(}
\def\){\right)}
\def\sls#1#2{(\f{#1}{#2})}
 
\def\Ls#1#2{\Big(\f{#1}{#2}\Big)}

\let \pro=\proclaim
\let \endpro=\endproclaim
\begin{document}
\par\q\par\q
 \centerline {\bf
On the properties of even and odd sequences}
$$\q$$
\centerline{Zhi-Hong Sun} $$\q$$ \centerline{School of Mathematical
Sciences, Huaiyin Normal University,} \centerline{Huaian, Jiangsu
223001, P.R. China} \centerline{E-mail: zhihongsun@yahoo.com}
\centerline{Homepage: http://www.hytc.edu.cn/xsjl/szh}

\abstract{\par\q In this paper we continue to investigate the
properties of those sequences $\{a_n\}$ satisfying the condition
$\sum_{k=0}^n\binom nk(-1)^ka_k=\pm a_n$ $(n\ge 0)$. As applications
we deduce new recurrence relations and congruences for Bernoulli and
Euler numbers.
\par\q
 }
\endabstract
\footnotetext[1] {The author is supported by the Natural Sciences
Foundation of China (grant 11371163).}
\section*{1. Introduction}
\par\q\ The classical binomial inversion formula states that
$a_n=\sum\limits_{k=0}^n \binom n k(-1)^kb_k$\ $(n=0,1,2,\ldots)$
 if and only if $b_n=\sum\limits_{k=0}^n \binom n k(-1)^k a_k
 (n=0,1,2,\ldots)$. Following [8] we continue to study those sequences $\{a_n\}$
with the property $\sum_{k=0}^n\b nk(-1)^ka_k=\pm a_n$
$(n=0,1,2,\ldots)$.
 \pro{Definition 1.1}  If a sequence $\{a_n\}$
satisfies the relation
$$\sum_{k=0}^n\b nk(-1)^ka_k=a_n\q (n=0,1,2,\ldots),$$
 we say that $\{a_n\}$ is an even sequence.  If $\{a_n\}$ satisfies the relation
$$\sum_{k=0}^n\b nk(-1)^ka_k=-a_n\q (n=0,1,2,\ldots),$$
 we say that $\{a_n\}$ is an odd sequence.\endpro
 \par From [8, Theorem 3.2] we know that $\{a_n\}$ is an even (odd)
sequence if and only if
$\t{e}^{-x/2}\sum_{n=0}^{\infty}a_n\f{x^n}{n!}$ is an even (odd)
function. Throughout this paper,  $S^+$ denotes the set of even
sequences, and  $S^-$ denotes the set of odd sequences.  In [8] the
author stated that
$$\Big\{\f 1{2^n}\Big\},\ \Big\{\b{n+2m-1}m^{-1}\Big\},
\ \Big\{\b{2n}n2^{-2n}\Big\}, \ \Big\{(-1)^n\int_0^{-1}\b
xndx\Big\}\in S^+.$$
\par Let $\{B_n\}$ be the Bernoulli numbers given by
$B_0=1$ and $\sum_{k=0}^{n-1}\b nkB_k=0$ $(n\ge 2)$. It is well
known that $B_1=-\f 12$ and $B_{2m+1}=0$ for $m\ge 1$. Thus,
$$\sum_{k=0}^n\b
nk(-1)^k\cdot (-1)^kB_k=B_n+ \sum_{k=0}^{n-1}\b nkB_k =(-1)^nB_n$$
and so $\{(-1)^nB_n\}\in S^+$ as claimed in [8].
\par  The Euler numbers
$\{E_n\}$ is
 defined by
$\f{2\t{e}^t}{\t{e}^{2t}+1}=\sum_{n=0}^{\infty}E_n\f{t^n}{n!}\
(|t|<\f{\pi}2)$, which is equivalent to (see [3]) $E_0=1,\
E_{2n-1}=0$ and $\sum_{r=0}^n\b{2n}{2r}E_{2r}=0\ (n\ge 1).$
  It is clear that
$$\align \sum_{n=0}^{\infty}\f{E_n-1}{2^n}\cdot\f{t^n}{n!}&
=\sum_{n=0}^{\infty}E_n\f{(t/2)^n}{n!}-\sum_{n=0}^{\infty}\f{(t/2)^n}{n!}
\\&=\f{2\t{e}^{\f t2}}{\t{e}^t+1}-\t{e}^{\f t2}
=\t{e}^{\f t2}\cdot \f{1-\t{e}^t}{1+\t{e}^t}\q\big(|t|<\pi\big).
\endalign$$
As $\f{1-\t{e}^{-t}}{1+\t{e}^{-t}}=\f{\t{e}^t-1}{\t{e}^t+1}$, we see
that $\t{e}^{-\f
t2}\sum_{n=0}^{\infty}\f{E_n-1}{2^n}\cdot\f{t^n}{n!}$ is an odd
function. Thus $\big\{\f{E_n-1}{2^n}\big\}$ is an odd sequence.
\par For two numbers $b$ and $c$, let $\{U_n(b,c)\}$ and $\{V_n(b,c)\}$ be
the Lucas sequences given by
$$U_0(b,c)=0,\ U_1(b,c)=1,\ U_{n+1}(b,c)=bU_n(b,c)-cU_{n-1}(b,c)\ (n\ge 1)$$
and
$$V_0(b,c)=2,\ V_1(b,c)=b,\ V_{n+1}(b,c)=bV_n(b,c)-cV_{n-1}(b,c)\ (n\ge 1).$$
 It is well
known that (see [12]) for $b^2-4c\not=0$,
$$U_n(b,c)= \f
1{\sqrt{b^2-4c}}\Big\{\Big(\f{b+\sqrt{b^2-4c}}2\Big)^n
-\Big(\f{b-\sqrt{b^2-4c}}2\Big)^n\Big\}$$ and
$$V_n(b,c)=\Ls{b+\sqrt{b^2-4c}}2^n+\Ls{b-\sqrt{b^2-4c}}2^n.$$
From this one can easily see that for $b(b^2-4c)\not=0$,
$\{U_n(b,c)/b^n\}$ is an odd sequence and $\{V_n(b,c)/b^n\}$ is an
even sequence.
\par Let $\{A_n\}$ be an even sequence or an odd sequence. In Section 2 we deduce
new recurrence formulas for $\{A_n\}$ and give a criterion for
polynomials $P_m(x)$ with the property $P_m(1-x)=(-1)^mP_m(x)$, in
Section 3 we establish a transformation formula for $\sum_{k=0}^n \b
nkA_k$,  in Section 4 we give a general congruence involving $A_n$
modulo $p^2$, where $p>3$ is a prime. As applications we establish
new recurrence formulas and congruences for Bernoulli and Euler
numbers. Here are some typical results:
\par $\star$ If $\{A_n\}$ is an even sequence and $n$ is odd, then
$$\sum_{k=0}^n\b{\f {n}2}k(-1)^kA_{n-k}=0\qtq{and}
\sum_{k=0}^n\b nk\b{n+k}k(-1)^kA_k=0.$$
\par $\star$ If $\{A_n\}$ is an odd
sequence, then $\sum_{k=0}^n\b nk(-1)^k A_{2n-k}=0$ for
$n=0,1,2,\ldots.$
\par$\star$ If $\{A_n\}$ is an even sequence and $\lambda$ is a real number, then
$$\sum_{k=0}^{2n+1}\b{2n-\lambda}{2n+1-k}\b{\lambda}k2^kA_k=0
\q(n=0,1,2,\ldots).$$
 \par $\star$ Let $m$ be a positive integer and $P_m(x)=\sum_{k=0}^m
a_kx^{m-k}$. Then
$$P_m(1-x)=(-1)^mP_m(x)\iff \sum_{k=0}^n\b nk\f{a_k}{\b mk}
=(-1)^n\f{a_n}{\b mn}\ (n=0,1,\ldots,m).$$
\par $\star$ Let $p$ be an odd prime, and let $\{A_k\}$ be an
odd sequence of rational $p$-integers. Then
$\sum_{k=1}^{p-1}\f{A_k}{p+k}\e 0\mod{p^2}.$
\par $\star$ Suppose that $\{a_n\}\in S^+$ with $a_0\not=0$ and $A_n=\f
1{(n+1)(n+2)}\sum_{k=0}^na_k$\q$ (n\ge 0)$. Then $\{A_n\}\in S^+$.
\par $\star$ Let $[x]$ be the greatest integer not exceeding $x$. For $n=3,4,5,\ldots$ we have
$$\sum_{r=1}^{[\f{n+1}2]}\b n{2r-1}(2n-2r+1)B_{2n-2r}=0\qtq{and}
\sum_{k=0}^{[n/2]}\b n{2k}2^{2k}E_{2n-2k}=(-1)^n.$$
\par In addition to the above notation throughout this paper we use the following notation:
$\Bbb N\f{\q}{\q}$the set of positive integers, $\Bbb
R\f{\q}{\q}$the set of real numbers, $\Bbb Z_p\f{\q}{\q}$the set of
those rational numbers whose denominator is coprime to $p$, $\sls
ap\f{\q}{\q}$the Legendre symbol.

\section*{2. Recurrence formulas for even and odd sequences}

 For $x,y\in\Bbb R$ and $n\in\{0,1,2,\ldots\}$ it is well known that
 $$\sum_{k=0}^n\b xk\b y{n-k}=\b{x+y}n.$$
 This is called Vandermonde's identity.
 Let $a_n=\nm b_{n-k}\quad (n=0,1,2,\ldots)$. Using Vandermonde's identity
 we see that
$$\aligned &\sum_{k=0}^n\binom {n-m}k(-1)^{n-k}
a_{n-k}
\\&=\nm\sum_{j=0}^{n-k}\b{n-k-m}{j}(-1)^{n-k-j}b_{n-k-j}
\\&=\sum_{s=0}^n\b{n-m}{n-s}(-1)^s\sum_{j=0}^s\b{s-m}{j}(-1)^{s-j}
b_{s-j}
\\&=\sum_{s=0}^n\b{n-m}{n-s}\sum_{r=0}^s\b{m-r-1}{s-r}b_r
\\&=\sum_{r=0}^n\sum_{s=r}^n\b{n-m}{n-s}\b{m-r-1}{s-r}b_r
\\&=\sum_{r=0}^n\b{n-r-1}{n-r}b_r
=b_n\qquad(n=0,1,2,\ldots).
\endaligned$$

 \proclaim{Lemma 2.1} Let
$m,p\in\Bbb R$ and $\sum_{k=0}^n\b{n-m}k(-1)^{n-k}a_{n-k}=\pm a_n$
for $n=0,1,2,\ldots.$ Then
$$\sum_{k=0}^n\b{n-p-m}k(-1)^{n-k}a_{n-k}
=\pm\sum_{k=0}^n\b pk(-1)^ka_{n-k}\qtq{for}n=0,1,2,\ldots.$$
\endproclaim
Proof. Using Vandermonde's identity we see that
$$\align \sum_{k=0}^n&\b{n-p-m}k(-1)^{n-k}a_{n-k}
\\&=\sum_{k=0}^n\b{n-p-m}{n-k}(-1)^ka_k
\\& =\pm\sum_{k=0}^n\b{n-p-m}{n-k}(-1)^k\sum_{r=0}^k\b{k-m}{k-r}(-1)^ra_r
\\&=\pm\sum_{r=0}^n\Big\{\sum_{k=r}^n\b{n-p-m}{n-k}(-1)^{k-r}\b{k-m}
{k-r}\Big\}a_r
\\&=\pm\sum_{r=0}^n\Big\{\sum_{k=r}^n\b{n-p-m}{n-k}\b{m-1-r}{k-r}\Big\}a_r
\\&=\pm\sum_{r=0}^n\Big\{\sum_{s=0}^{n-r}\b{n-p-m}{n-r-s}\b{m-1-r}s\Big\}a_r
\\&=\pm\sum_{r=0}^n\b{n-p-r-1}{n-r}a_r
=\pm\sum_{r=0}^n\b{p}{n-r}(-1)^{n-r}a_r
\\&=\pm\sum_{k=0}^n\b pk(-1)^ka_{n-k}.
\endalign$$
So the lemma is proved.
 \pro{Theorem 2.1} Let $\{A_n\}$ be an even sequence.  For $n=1,3,5,\ldots$ we have
$$\sum_{k=0}^n\b{\f {n}2}k(-1)^kA_{n-k}=0.$$\endpro
Proof. Putting $m=0$, $p=n/2$ and $a_n=A_n$ in Lemma 2.1 we deduce
the result.
\pro{Corollary 2.1} For $n=1,3,5,\ldots$ we have
$$\sum_{k=0}^n\b{\f n2}kB_{n-k}=0\q\t{and}\q
\sum_{k=0}^n\b{\f n2}k\f{(2^{n-k+1}-1)B_{n-k+1}}{n-k+1}=0.$$
\endpro
Proof. From [8] we know that $\{(-1)^nB_n\}\in S^+$ and
$\{(-1)^{n+1}(2^{n+1}-1)B_{n+1}/(n+1)\}\in S^+$. Thus the result
follows from Theorem 2.1.
 \pro{Lemma 2.2} If $\{a_n\}\in S^+$, then
 $\{na_{n-1}\},\{\f{a_{n+1}}{n+1}\}\in S^-$.
If $\{a_n\}\in S^-$, then
 $\{na_{n-1}\},\{\f{a_{n+1}}{n+1}\}\in S^+$.
 \endpro
 Proof. Suppose that $\sum_{k=0}^n\b nk(-1)^ka_k=\pm a_n$ for
  $n=0,1,2,\ldots.$ Since
$$\align \sum_{k=0}^n\b nk(-1)^kka_{k-1}&
=\sum_{k=1}^nn\b{n-1}{k-1}(-1)^ka_{k-1}
\\&=-n\sum_{r=0}^{n-1}\b{n-1}r(-1)^ra_r=\mp na_{n-1}\endalign$$
and $$\align\sum_{k=0}^n\b nk(-1)^k\f{a_{k+1}}{k+1}&=\f
1{n+1}\sum_{k=0}^n\b{n+1}{k+1}(-1)^ka_{k+1}
\\&=-\f 1{n+1}\sum_{r=0}^{n+1}\b{n+1}r(-1)^ra_r=\mp\f{a_{n+1}}{n+1},
\endalign$$ we see that the result is true.
 \pro{Theorem 2.2} If $\{A_n\}$ is an odd
sequence, then
$$\sum_{k=0}^n\b nk(-1)^k A_{2n-k}=0\qtq{for}n=0,1,2,\ldots.$$
If $\{A_n\}$ is an even sequence, for $n=0,1,2,\ldots$ we have
$$\sum_{k=0}^n\b nk(-1)^k(2n-k)A_{2n-k-1}=0\qtq{and}
\sum_{k=0}^n\b nk(-1)^k\f{A_{2n-k+1}}{2n-k+1}=0 .$$
\endpro
Proof. We first assume that$\{A_n\}$ is an odd sequence. Putting
$m=0$, $p=n/2$ and $a_n=A_n$ in Lemma 2.1 we see that
$$\sum_{k=0}^n\b{\f n2}k(-1)^{n-k}A_{n-k}=-\sum_{k=0}^n\b{\f
n2}k(-1)^kA_{n-k}.$$ Thus, for even $n$ we have
$$\sum_{k=0}^{n/2}\b{\f n2}k(-1)^kA_{n-k}
=\sum_{k=0}^n\b{\f n2}k(-1)^kA_{n-k}=0.$$ Replacing $n$ with $2n$ we
get $\sum_{k=0}^n\b nk(-1)^kA_{2n-k}=0.$
\par Now we assume that $\{A_n\}\in S^+$. By Lemma 2.2,
$\{nA_{n-1}\},\{\f{A_{n+1}}{n+1}\}\in S^-$. Thus applying the above
we deduce the remaining result.
\par\q

\pro{Corollary 2.2} For $n=0,1,2,\ldots$ we have
$$\sum_{k=0}^{[n/2]}\b n{2k}2^{2k}E_{2n-2k}=(-1)^n.$$
\endpro
Proof. As $\{\f{E_n-1}{2^n}\}\in S^-$, taking $A_n=\f{E_n-1}{2^n}$
in Theorem 2.2 we obtain $$\sum_{k=0}^n\b
nk(-1)^k\f{E_{2n-k}-1}{2^{2n-k}}=0.$$ That is,
$$\sum_{k=0}^n\b nk(-2)^kE_{2n-k}=\sum_{k=0}^n\b nk(-2)^k=(-1)^n.$$
To see the result, we note that $E_{2m-1}=0$ for $m\ge 1$.
\pro{Corollary 2.3} For $n=3,4,5,\ldots$ we have
$$\sum_{r=1}^{[\f{n+1}2]}\b n{2r-1}(2n-2r+1)B_{2n-2r}=0.$$
\endpro
Proof. As $\{(-1)^nB_n\}\in S^+$, taking $A_n=(-1)^nB_n$ in Theorem
2.2 we see that $\sum_{k=0}^n\b nk (2n-k)B_{2n-k-1}=0$. To see the
result, we note that $B_{2m+1}=0$ for $m\ge 1$.

 \pro{Corollary 2.4} For $n=2,3,4,\ldots$ we have
 $$\sum_{r=0}^{[n/2]}\b n{2r}(2^{2n-2r}-1)B_{2n-2r}=0.$$
 \endpro
 Proof. The Bernoulli polynomials $\{B_n(x)\}$ are given
 by  $B_n(x)= \sum_{k=0}^n\b nkB_kx^{n-k}$. It is well known
 ([3]) that $B_n(1-x)=(-1)^nB_n(x)$. From [3, p.248] we also have
$B_{2k}(\f 12)=(2^{1-2k}-1)B_{2k}$. Thus $B_{2k+1}(\f 12)=0$ and so
$B_n(\f 12)=(2^{1-n}-1)B_n$. Hence,
$$\align \sum_{k=0}^n\b nk(2^k-1)B_k
&=2^n\sum_{k=0}^n\b nkB_k\Ls 12^{n-k}-\sum_{k=0}^n\b nkB_k
\\&=2^nB_n\Ls 12-(-1)^nB_n=(2-2^n)B_n-(-1)^nB_n
\\&=(-1)^n(1-2^n)B_n.
\endalign$$
That is, $\{(-1)^n(2^n-1)B_n\}$ is an odd sequence. Now applying
Theorem  2.2 we obtain $\sum_{k=0}^n\b nk(2^{2n-k}-1)B_{2n-k}=0$. To
see the result, we note that $B_{2m+1}=0$ for $m\ge 1$.

 \pro{Lemma 2.3} Let $m,n\in\Bbb N$ with
$m\le n$. Then$$\sum_{k=m}^n\b nk\b{n+k}k(-1)^{n-k}\b km=\b
nm\b{m+n}m.$$
\endpro
Proof. Using Vandermonde's identity we see that
$$\sum_{k=0}^n\b{n-m}{n-k}\b{-n-1}k=\b{-m-1}n.$$
Since $\b xk=(-1)^k\b{x+k-1}k$, we have
$$\sum_{k=0}^n\b{n-m}{n-k}(-1)^k\b{n+k}k=(-1)^n\b{m+n}n.$$
That is,
$$\sum_{k=m}^n\b{n-m}{k-m}(-1)^{n-k}\b{n+k}k=\b{m+n}n.$$
Hence
$$\sum_{k=m}^n\b nm\b{n-m}{k-m}(-1)^{n-k}\b{n+k}k=\b nm\b{m+n}m.$$
 As $\b nk\b km=\b nm\b{n-m}{k-m}$, we obtain the result.

 \pro{Lemma 2.4} Let
$\{a_k\}$ be a given sequence. For $n\in\Bbb N$ we have
$$\sum_{k=0}^n\b nk\b{n+k}k\Big(a_k-(-1)^{n-k}\sum_{s=0}^k
\b ksa_s\Big)=0.$$

\endpro
Proof. Using Lemma 2.3 we see that
$$\align &\sum_{k=0}^n\b nk\b{n+k}k\Big(a_k-(-1)^{n-k}\sum_{s=0}^k\b ksa_s\Big)
\\&=\sum_{m=0}^na_m\Big(\b nm\b{n+m}m-
\sum_{k=m}^n\b nk\b{n+k}k(-1)^{n-k}\b km\Big)
\\&=\sum_{m=0}^na_m\cdot 0=0.\endalign$$
\pro{Theorem 2.3} If $\{A_n\}$ is an even sequence and $n$ is odd,
or if   $\{A_n\}$ is an odd sequence and $n$ is even, then
$$\sum_{k=0}^n\b nk\b{n+k}k(-1)^kA_k=0.$$
\endpro
Proof. Putting $a_k=(-1)^kA_k$ in Lemma 2.4 we obtain the result.
\proclaim{Theorem 2.4} Suppose that $m$ is a nonnegative integer.
Then $$\sum_{k=0}^n\b{n-m-1}k(-1)^{n-k} a_{n-k} =\pm a_n\
(n=0,1,2,\ldots)$$ if and only if
$$\sum_{k=0}^n\b nk\frac{a_k}{\b mk}
=\pm(-1)^n\frac{a_n}{\b mn}\ (n=0,1,\ldots,m)$$ and
$$\n a_{k+m+1}=\pm(-1)^{m+1}a_{n+m+1}\ (n=0,1,2,\ldots).$$\endproclaim
Proof. For $n=0,1,\ldots,m$ we have $\b mn\neq 0$. Set
$A_n=(-1)^n\frac{a_n}{\b mn}$. As $\b{n-m-1}k\b m{n-k}=(-1)^k\b nk\b
mn$, we see that
$$\align &\sum_{k=0}^n\b{n-m-1}k(-1)^{n-k}a_{n-k}
\\&=\sum_{k=0}^n\b{n-m-1}k\b m{n-k}A_{n-k}
=\b mn\sum_{k=0}^n\b nk(-1)^kA_{n-k} \\&=(-1)^n\b mn\sum_{k=0}^n\b
nk(-1)^kA_k.\endalign$$ Thus,
 $$\sum_{k=0}^n\b{n-m-1}k(-1)^{n-k}a_{n-k}
 =\pm a_n\lra \sum_{k=0}^n\b nk(-1)^kA_k=\pm A_n.\tag 2.1$$
This together with the fact that
$$\align&\sum_{k=0}^{n+m+1}\b{n+m+1-m-1}k(-1)^{n+m+1-k}a_{n+m+1-k}
\\&\ =\sum_{k=0}^n\b nk(-1)^{n-k+m+1}a_{n-k+m+1}
\\&\ =\sum_{r=0}^n\b nr(-1)^{r+m+1}a_{r+m+1}\qquad(n=0,1,2,\ldots)\endalign$$
yields the result.
 \par For any sequence $\an$ the formal power
series $\sum_{n=0}^{\infty}a_nx^n$ is
 called the generating function of $\an$.
 \proclaim{Lemma 2.5} Let
 $\{a_n\}$ be a given sequence, $a(x)=\sum_{n=0}^{\infty}a_nx^n$ and
 $m\in\Bbb R$.
Then
$$\align &(1-x)^ma\Big(\f x{x-1}\Big)=\pm a(x)
\\&\iff\sum_{k=0}^n\b{n-m-1}k(-1)^{n-k}a_{n-k}=\pm a_n\ (\ncdots).
\endalign$$
 \endproclaim
   Proof. Clearly,
 $$\align  (1-x)^ma\Big(\f x{x-1}\Big)&
 =\sum_{r=0}^{\infty}(-1)^ra_rx^r(1-x)^{m-r}
  \\&=\sum_{r=0}^{\infty}(-1)^ra_rx^r\sum_{k=0}^{\infty}
  \b {m-r}k(-x)^k
 \\&=\sum_{n=0}^{\infty}\bigg(\sum_{k=0}^n(-1)^{n-k}a_{n-k}
 \b{m-(n-k)}k
 (-1)^k\bigg)x^n
 \\&=\sum_{n=0}^{\infty}\bigg(\sum_{k=0}^n\b{n-m-1}k(-1)^{n-k}
 a_{n-k}\bigg) x^n.\endalign$$
 Thus the result follows.

 \proclaim{Theorem 2.5} Let $m\in\Bbb N$,
  $P_m(x)=\sum_{k=0}^m
a_kx^{m-k}$ and $P_m^*(x)=\sum_{k=0}^ma_kx^k$. Then the following
statements are equivalent:
\par \abc a $(1-x)^mP_m^*\big(\frac
x{x-1}\big)=\pm P_m^*(x)$.\par \abc b $P_m(1-x)=\pm(-1)^mP_m(x)$.
\par \abc c For $n=0,1,\ldots,m$ we have
$\sum_{k=0}^n\b{n-m-1}k(-1)^{n-k} a_{n-k}=\pm a_n.$
\par \abc d Set $a_n=0$ for $n>m$. Then
$\sum_{k=0}^n\b{n-m-1}k(-1)^{n-k} a_{n-k}=\pm a_n\
(n=0,1,2,\ldots).$
\par \abc e For $n=0,1,\ldots,m$
we have
$$\sum_{k=0}^n\b nk\frac{a_k}{\b mk}=\pm(-1)^n
\frac{a_n}{\b mn}.$$\endproclaim Proof. Since $P_m^*(x)=x^mP_m(\frac
1x)$ we see that
$$\align(1-x)^mP_m^*\Big(\frac x{x-1}\Big)=
\pm P_m^*(x)&\lra (-x)^mP_m\Big(1-\frac 1x\Big) =\pm
x^mP_m\Big(\frac 1x\Big)
\\&\lra\frac 1{x^m}\left((-1)^mP_m(1-x)\mp P_m(x)\right)=0
\\&\lra P_m(1-x)=\pm
(-1)^mP_m(x).\endalign$$ So (a) and (b) are equivalent. By Lemma
2.5, (a) is equivalent to (d).  Assume $a_{n+m+1}=0$ for $n\ge 0$.
Then
$$\align a_{n+m+1}&=0=\sum_{k=0}^n\b nk(-1)^{n-k+m+1}a_{m+1+n-k}
\\&=\sum_{k=0}^{m+n+1}\b{m+n+1-m-1}k (-1)^{m+n+1-k}a_{m+n+1-k}.\endalign$$
So (c) is equivalent to (d).  To complete the proof, we note that
(d) is equivalent to (e)
 by Theorem 2.4.

  \proclaim{Corollary 2.5} Let $m\in\Bbb R$. Then $$\sum_{k=0}^n\b{n-m-1}k(-1)^{n-k}
a_{n-k}=\pm a_n\ (\ncdots)$$ if and only if
$$\sum_{k=0}^p\b pk(-1)^k\frac{\b {n-m-1}{n-k}}{\b nk}a_k=\pm\frac{\b{n-m-1}{n-p}}
{\b np}a_p\tag 2.2$$ for every nonnegative integer $n$ and every
$p\in\{0,1,\ldots,n\}.$
\endproclaim
Proof. If (2.2) holds, taking $p=n$ in (2.2) we see that
$$\sum_{k=0}^n\b{n-m-1}k(-1)^{n-k} a_{n-k}=\pm a_n\qq(\ncdots).\tag 2.3 $$
 Set
$P_n(x)=\sum_{k=0}^n\b{n-m-1}{n-k}(-1)^ka_kx^{n-k}$. Then
$$\align &P_n(1-x)\\&=\sum_{k=0}^n\b{n-m-1}k(-1)^{n-k} a_{n-k}\sum_{r=0}^k\b kr (-x)^r
\\&=\sum_{r=0}^n\b {n-m-1}r(-1)^r\sum_{k=r}^n\b {n-m-1-r}{k-r}(-1)^{n-k}a_{n-k}
x^r
\\&=(-1)^n\sum_{r=0}^n\b{n-m-1}r(-1)^{n-r}
\left(\sum_{s=0}^{n-r}\b {n-r-m-1}s(-1)^{n-r-s}a_{n-r-s}\right)
x^r.\endalign $$ Thus,
$$\aligned &P_n(1-x)=\pm (-1)^nP_n(x)
\\&\Longleftrightarrow\sum_{s=0}^k\b{k-m-1}s
(-1)^{k-s}a_{k-s}=\pm a_k\ \t{for}\  k=0,1,\ldots,n.
\endaligned \tag 2.4$$
If (2.3) holds, from (2.4) we see that $P_n(1-x)=\pm(-1)^nP_n(x)$
and so (2.2) holds by Theorem 2.5. This completes the proof.

\proclaim{Theorem 2.6} Suppose that $m,p\in\Bbb R$, $\sum_{k=0}^n
\b{n-m}k(-1)^{n-k}a_{n-k}=(-1)^{\a}a_n$ and $\sum_{k=0}^n
\b{n-p}k(-1)^{n-k}b_{n-k}=(-1)^{\beta} b_n\q(\ncdots)$. If $n$ is a
nonnegative integer such that $\a+\beta+n$ is odd, then
$$\sum_{k=0}^n\frac{\b
{n-m}k\b{n-p}{n-k}}{\b nk}(-1)^k a_{n-k}b_k =0.$$
\endproclaim
Proof.
 Set $T_n=\sum_{k=0}^n\frac{\b {n-m}k\b{n-p}{n-k}}{\b nk}(-1)^k
a_{n-k}b_k.$
  From Corollary 2.5 we know that
$$\sum_{r=0}^k\b kr(-1)^r\frac{\b {n-p}{n-r}}{\b nr}b_r=(-1)^{\beta}\frac
{\b{n-m}{n-k}}{\b nk}b_k$$ for every nonnegative integer $n$ and
$k\in\{0,1,\ldots,n\}$. Thus,
$$\align T_n=&\sum_{k=0}^n\b{n-m}k(-1)^ka_{n-k}\cdot(-1)^{\beta}\sum_{r=0}
^k\b kr(-1)^r\frac{\b{n-p}{n-r}}{\b nr}b_r
\\=&(-1)^{\beta}\sum_{r=0}^n\frac{\b{n-p}{n-r}}{\b nr}b_r\sum_{k=r}^n\b
{n-m}k\b kr(-1)^{k-r}a_{n-k}
\\=&(-1)^{\beta}\sum_{r=0}^n\frac{\b{n-p}{n-r}}{\b nr}b_r\b{n-m}r
\sum_{k=r}^n\b{n-m-r}{k-r}(-1)^{k-r}a_{n-k}
\\=&(-1)^\beta\sum_{r=0}^n\frac{\b{n-m}r\b{n-p}r}{\b nr}b_r\sum_{s=0}
^{n-r}\b{n-r-m}s(-1)^sa_{n-r-s}
\\=&(-1)^{\a+\beta+n}\sum_{r=0}^n\frac{\b{n-m}r\b{n-p}r}{\b nr}(-1)^r
a_{n-r}b_r
\\=&(-1)^{\a+\beta+n}T_n\qq(\ncdots).\endalign$$
Hence $T_n=0$ when $\a+\beta+n$ is odd. This completes the proof.

\pro{Theorem 2.7} Let $\lambda$ be a real number.
\par $(\hbox{\rm i})$ If $\{a_n\},\{A_n\}\in
S^+$ or $\{a_n\},\{A_n\}\in S^-$, then
$$\sum_{k=0}^{2n+1}\b{2n-\lambda}{2n+1-k}\b{\lambda}ka_{2n+1-k}A_k=0
\q(n=0,1,2,\ldots).$$
\par $(\hbox{\rm ii})$ If $\{a_n\}\in S^+$ and $\{A_n\}\in
S^-$, then
$$\sum_{k=0}^{2n}\b{2n-1-\lambda}{2n-k}\b{\lambda}ka_{2n-k}A_k=0
\q(n=0,1,2,\ldots).$$
\endpro
Proof. Suppose $\sum_{k=0}^n\b nk(-1)^ka_k=(-1)^{\a}a_n$ and
$\sum_{k=0}^n\b nk(-1)^kA_k=(-1)^{\beta}A_n$ for $n=0,1,2,\ldots$.
Set $b_n=(-1)^n\b{\lambda}nA_n$. By (2.1) we have
$$\sum_{k=0}^n\b{n-1-\lambda}k(-1)^{n-k}b_{n-k}=(-1)^{\beta}b_n.$$ Now taking $m=0$ and
$p=\lambda+1$ in Theorem 2.6 we see that if $2\nmid\a+\beta+n$, then
$$\sum_{k=0}^n\b{n-1-\lambda}{n-k}a_{n-k}\b{\lambda}kA_k
=\sum_{k=0}^n\b{n-1-\lambda}{n-k}(-1)^ka_{n-k}b_k=0.$$ This yields
the result.
 \pro{Corollary 2.6} Let $\lambda$ be a real number.
\par $(\hbox{\rm i})$ If $\{A_n\}\in S^+$, then
$$\sum_{k=0}^{2n+1}\b{2n+1-\lambda}{2n+2-k}\b{\lambda}kA_k=0\q(n=0,1,2,\ldots).$$
\par $(\hbox{\rm ii})$ If $\{A_n\}\in S^-$, then
$$\sum_{k=0}^{2n}\b{2n-\lambda}{2n+1-k}\b{\lambda}kA_k=0\q(n=0,1,2,\ldots).$$
\endpro
Proof. As $\{\f 1{n+1}\}\in S^+$, putting $a_n=\f 1{n+1}$ in Theorem
2.7 we deduce the result.
 \pro{Corollary 2.7} Let $\lambda$
be a real number.
\par $(\hbox{\rm i})$ If $\{A_n\}\in S^+$, then
$$\sum_{k=0}^{2n+1}\b{2n-\lambda}{2n+1-k}\b{\lambda}k2^kA_k=0
\q(n=0,1,2,\ldots).$$
\par $(\hbox{\rm ii})$ If $\{A_n\}\in S^-$, then
$$\sum_{k=0}^{2n}\b{2n-1-\lambda}{2n-k}\b{\lambda}k2^kA_k=0
\q(n=0,1,2,\ldots).$$
\endpro
Proof. As $\{\f 1{2^n}\}\in S^+$, putting $a_n=\f 1{2^n}$ in Theorem
2.7 we deduce the result.

\proclaim{Theorem 2.8} Let $\{A_n\}\in S^+$ with
$A_0=\ldots=A_{l-1}=0$ and $A_l\not=0\ (l\ge 1)$. Then
$$\Big\{\f{A_{n+l}}{(n+1)(n+2)\cdots(n+l)}\Big\}\in S^+.$$
\endproclaim
Proof. Assume $a_n=A_{n+l}$. Let $a(x)$ and $A(x)$ be the generating
functions of $\{a_n\}$ and $\{A_n\}$ respectively. Then clearly
$A(x)=x^la(x)$. Since $A_l=\sum_{k=0}^l\b lk(-1)^kA_k=(-1)^lA_l$ we
see that $2\mid l$. Thus, applying Lemma 2.5 and (2.1) we see that
$$\align\{A_n\}\in S^+&\Leftrightarrow A(\f x{x-1})=(1-x)A(x)
\Leftrightarrow a\Big(\f x{x-1}\Big)=(1-x)^{l+1}a(x)
\\&\Leftrightarrow \sum_{r=0}^n\b{n+l}r(-1)^{n-r}a_{n-r}=a_n
\ (\ncdots)
\\&\Leftrightarrow\Big\{\f{(-1)^na_n}{\b{-l-1}n}\Big\}\in S^+.\endalign$$
 Note that
$$(-1)^n\f{a_n}{\b{-l-1}n}=\f {a_n}{\b{n+l}l}=\f{A_{n+l}}{(n+1)(n+2)\cdots(n+l)}
\cdot l!.$$ We then obtain the result. \proclaim{Theorem 2.9}
Suppose that $\{a_n\}\in S^+$ with $a_0\not=0$ and $A_n=\f
1{(n+1)(n+2)}\sum_{k=0}^na_k$\q$ (n\ge 0)$. Then $\{A_n\}\in S^+$.
\endproclaim
Proof. Let $b_0=b_1=0$ and $b_n=\sum_{k=0}^{n-2}a_k\ (n\ge 2)$. Let
$a(x)=\sum_{n=0}^{\infty}a_nx^n$ and
$b(x)=\sum_{n=0}^{\infty}b_nx^n$.
 Then
$$\align b(x)&=\sum_{n=2}^{\infty}\Big(\sum_{k=0}^{n-2}a_k\Big)x^n
=x^2\sum_{n=0}^{\infty}\Big(\sum_{k=0}^na_k\Big)x^n
\\&=\f {x^2}{1-x}\sum_{n=0}^{\infty}a_nx^n=\f {x^2}{1-x}a(x).
\endalign$$
Thus, by [8, Theorem 3.1] or Lemma 2.5 (with $m=-1$) we have
$$\align \{a_n\}\in S^+&\iff (1-x)a(x)=a\Big(\f x{x-1}\Big)
 \iff (1-x)b(x)=b\Big(\f x{x-1}\Big)
 \\& \iff  \{b_n\}\in S^+.\endalign$$
 Since
$b_0=b_1=0$ and $b_2=a_0\not=0$, applying Theorem 2.8 we find that
$\{\f{b_{n+2}}{(n+1)(n+2)}\}\in S^+$. That is $\{A_n\}\in S^+$.

\pro{Theorem 2.10} Let $F$ be a given function. If $\{A_n\}$ is an
even sequence, then
$$\sum_{k=0}^n\b nk(-1)^kA_k\Big(\sum_{s=0}^k\b ks(-1)^s(F(s)-
F(n-s))\Big)=0\q(n=0,1,2,\ldots).$$ If $\{A_n\}$ is an odd sequence,
then
$$\sum_{k=0}^n\b nk(-1)^kA_k\Big(\sum_{s=0}^k\b ks(-1)^s(F(s)+
F(n-s))\Big)=0\q(n=0,1,2,\ldots).$$
\endpro
Proof. Suppose that $\sum_{k=0}^n\b nk(-1)^kA_k=\pm A_n$. From [7,
Lemma 2.1] we have
$$\sum_{k=0}^n\b nk(-1)^kf(k)A_k=\pm\sum_{k=0}^n\b
nk\Big(\sum_{r=0}^k\b kr(-1)^rF(n-k+r)\Big)A_k,$$ where
$f(k)=\sum_{s=0}^k\b ks(-1)^sF(s)$. Thus
$$\sum_{k=0}^n\b nk(-1)^kA_k\Big(f(k)\mp \sum_{s=0}^k\b
ks(-1)^sF(n-s)\Big)=0.$$ This yields the result. \pro{Corollary 2.8}
If $\{A_n\}\in S^+$, then
$$\sum_{k=0}^n \b nk(-1)^kA_k(1+x)^k(1-(-1)^nx^{n-k})=0\qtq{for}n=0,1,2,
\ldots.$$ If $\{A_n\}\in S^-$, then
$$\sum_{k=0}^n \b nk(-1)^kA_k(1+x)^k(1+(-1)^nx^{n-k})=0\qtq{for}n=0,1,2,
\ldots.$$
\endpro
Proof. Taking $F(s)=(-x)^s$ in Theorem 2.10 and then applying the
binomial theorem we obtain the result.
\par\q
\par From [7, (2.5)] we know that
$$\sum_{k=0}^n\b nk(-1)^kf(m+k)=\sum_{k=0}^m\b mk(-1)^kF(n+k),$$
where $F(r)=\sum_{s=0}^r\b rs(-1)^sf(s)$. Hence we have:
 \pro{Theorem 2.11} If $\{A_n\}$ is an even sequence, then
  for any nonnegative integers $m$ and $n$ we
have
$$\sum_{k=0}^n\b nk(-1)^kA_{k+m}=\sum_{k=0}^m\b mk(-1)^kA_{k+n}.$$
If $\{A_n\}$ is an odd sequence, then
  for any nonnegative integers $m$ and $n$ we
have
$$\sum_{k=0}^n\b nk(-1)^kA_{k+m}=-\sum_{k=0}^m\b mk(-1)^kA_{k+n}.$$
\endpro
From [8] we know that $\{1/\b{n+2r-1}r\}\in S^+$ for $r=1,2,\ldots$.
Thus, by Theorem 2.11 we have
$$\sum_{k=0}^n\b nk(-1)^k\f 1{\b{k+m+2r-1}r}
=\sum_{k=0}^m\b mk(-1)^k\f 1{\b{k+n+2r-1}r}.\tag 2.5$$
 Since $\{\f{U_n(b,c)}{b^n}\}\in S^-$ and $\{\f{V_n(b,c)}{b^n}\}\in
S^+$ for $b(b^2-4c)\not=0$, by Theorem 2.11 we have
$$\align&\sum_{k=0}^n\b nk(-1)^k\f{U_{k+m}(b,c)}{b^{k+m}}=-\sum_{k=0}^m\b mk
(-1)^k\f{U_{k+n}(b,c)}{b^{k+n}},\tag 2.6\\& \sum_{k=0}^n\b
nk(-1)^k\f{V_{k+m}(b,c)}{b^{k+m}} =\sum_{k=0}^m\b
mk(-1)^k\f{V_{k+n}(b,c)}{b^{k+n}} .\tag 2.7\endalign$$

\section*{3. A transformation formula for $\sum_{k=0}^n\b nkA_n$}

\pro{Lemma 3.1 ([8, Theorems 4.1 and 4.2])} Let $f$ be a given
function and $n\in\Bbb N$.
\par $(\t{\rm i})$ If $\{A_n\}$ is an even sequence, then
$$\sum_{k=0}^n\b
nk\Big(f(k)-(-1)^{n-k}\sum_{s=0}^k\b ksf(s)\Big)A_{n-k}=0.$$
\par $(\t{\rm ii})$ If $\{A_n\}$ is an odd sequence, then
$$\sum_{k=0}^n\b
nk\Big(f(k)+(-1)^{n-k}\sum_{s=0}^k\b ksf(s)\Big)A_{n-k}=0.$$
\endpro
\par We remark that a simple proof of Lemma 3.1 was given by
Wang[11].
 \pro{Theorem 3.1} Let $n\in\Bbb N$.
 If $\{A_m\}$ is an even sequence and $n$ is odd, or if
$\{A_m\}$ is an odd sequence and $n$ is even, then
$$\sum\Sb k=0\\3\mid k\endSb^n\b nkA_{n-k}
=\sum\Sb k=0\\3\mid n-k\endSb^n\b nkA_k=\f 13\sum_{k=0}^n\b nkA_k.$$
\endpro
Proof. Set $\omega=(-1+\sqrt{-3})/2$. If $\{A_m\}$ is an even
sequence and $n$ is odd, putting $f(k)=\o^k$ in Lemma 3.1 we obtain
$$\sum_{k=0}^n\b nk(\o^k-(-1)^{n-k}(1+\o)^k)A_{n-k}=0.$$
As $1+\o=-\o^2$, we  have $\sum_{k=0}^n\b
nk(\o^k+\o^{2k})A_{n-k}=0.$ Therefore,
$$\aligned &3\sum\Sb k=0\\3\mid k\endSb^n\b nkA_{n-k}-\sum_{k=0}^n\b
nkA_k
\\&=\sum_{k=0}^n\b
nk(1+\o^k+\o^{2k})A_{n-k}-\sum_{k=0}^n\b nkA_{n-k}
\\&=\sum_{k=0}^n\b
nk(\o^k+\o^{2k})A_{n-k}=0.\endaligned$$ The remaining part can be
proved similarly.
\endpro
\pro{Corollary 3.1 (Ramanujan [1,5])} For $n=3,5,7,\ldots$ we have
$$\sum\Sb k=0\\6\mid k-3\endSb^n\b nkB_{n-k}
=\cases -\f n6&\t{if $n\e 1\mod 6$,}
\\\f n3&\t{if $n\e 3,5\mod 6$.}
\endcases$$
\endpro
Proof. As $\{(-1)^nB_n\}\in S^+$, taking $A_n=(-1)^nB_n$ in Theorem
3.1 we obtain
$$\aligned\sum\Sb k=0\\3\mid k\endSb^n\b nk(-1)^{n-k}B_{n-k}
&=\f 13\sum_{k=0}^n\b nk(-1)^kB_k =\f 13\Big(\sum_{k=0}^n\b
nkB_k+n\Big)\\&=\f 13(n+B_n)=\f n3.\endaligned$$ To see the result,
we note that
$$\sum\Sb k=0\\3\mid k\endSb^n\b nk(-1)^{n-k}B_{n-k}-\sum\Sb k=0\\6\mid k-3\endSb^n\b nkB_k=
\cases -nB_1=\f n2&\t{if $n\e 1\mod 6$,}\\0&\t{if $n\e 3,5\mod 6$.}
\endcases$$

\pro{Corollary 3.2 (Ramanujan [5])} For $n=6,8,10,\ldots$ we have
$$\f 43(2^n-1)B_n+\sum_{k=1}^{[n/6]}\b n{6k}(2^{n-6k}-1)B_{n-6k}=\cases
-\f n6&\t{if $n\e 4\mod 6$,}
\\\f n3&\t{if $n\e 0,2\mod 6$.}
\endcases$$
\endpro
Proof. Since $\{(-1)^n(2^n-1)B_n\}$ is an odd sequence, by Theorem
3.1 we have
$$\aligned&\sum\Sb k=0\\3\mid k\endSb^n\b nk(-1)^{n-k}(2^{n-k}-1)B_{n-k}
\\&=\f 13\sum_{k=0}^n\b nk(-1)^k(2^k-1)B_k =\f 13\Big(\sum_{k=0}^n\b
nk(2^k-1)B_k+n\Big)\\&=\f 13(-(-1)^n(2^n-1)B_n+n)=\f n3-\f
13(2^n-1)B_n.\endaligned$$ On the other hand,
$$\sum\Sb k=0\\3\mid k\endSb^n\b nk(-1)^{n-k}(2^{n-k}-1)B_{n-k}
-\sum\Sb k=0\\6\mid k\endSb^n\b nk(2^{n-k}-1)B_{n-k} =\cases
-nB_1=\f n2&\t{if $6\mid n-4$,}
\\0&\t{if $6\nmid n-4$.}
\endcases$$
Thus the result follows.

\pro{Corollary 3.3 (Lehmer [2])} For $n=6,8,10,\ldots$ we have
$$E_n+3\sum_{k=1}^{[n/6]}\b
n{6k}2^{6k-2}E_{n-6k}=\f{1+(-3)^{n/2}}2.$$
\endpro
Proof. Since $\{(E_n-1)/2^n\}$ is an odd sequence, by Theorem 3.1
and the fact $E_{2k+1}=0$ we have
$$\aligned\sum\Sb k=0\\6\mid k\endSb^n\b nk\f{E_{n-k}}{2^{n-k}}
-\sum\Sb k=0\\3\mid k\endSb^n\b nk\f 1{2^{n-k}} &=\sum\Sb k=0\\3\mid
k\endSb^n\b nk\f{E_{n-k}-1}{2^{n-k}} =\f 13\sum_{k=0}^n\b
nk\f{E_k-1}{2^k}
\\&=\f 13\Big\{\sum_{k=0}^n\b
nk(-1)^k\f{E_k-1}{2^k}+\Big(1-\f 12\Big)^n-\Big(1+\f
12\Big)^n\Big\}\\&=\f
13\Big\{-\f{E_n-1}{2^n}+\f{1-3^n}{2^n}\Big\}=\f{2-3^n-E_n}{3\cdot
2^n}.\endaligned$$ Observe that
$$\aligned\sum\Sb k=0\\3\mid k\endSb^n\b nk2^k
&=\sum_{k=0}^n\b nk2^k\cdot \f 13(1+\omega^k+\omega^{2k})
\\& =\f 13\big((1+2)^n+(1+2\omega)^n+(1+2\omega^2)^n\big)
\\&=\f 13\big(3^n+(\sqrt{-3})^n+(-\sqrt{-3})^n\big)=\f 13(3^n+2\cdot
(-3)^{\f n2}).
\endaligned$$
 We obtain
$$\f 43E_n+\sum\Sb k=1\\6\mid k\endSb^n\b nk2^kE_{n-k}
=\f{2-3^n}3+\f{3^n+2\cdot (-3)^{\f n2}}3=\f 23\big(1+(-3)^{\f
n2}\big).$$ This yields the result.
\par{\bf Remark 3.1} Compared with known
proofs of Corollaries 3.1-3.3 (see [1,2,5]), our proofs are simple
and natural.
 \par Now we introduce new sequence $\{S_n\}$ defined
by
$$S_n+\sum_{k=0}^n\b nkS_k=2\q(n=0,1,2,\ldots).\tag 3.1$$
 The first few values of $S_n$ are shown below:
$$\align &S_0=1,\ S_1=\f 12,\ S_3=-\f 14,\ S_5=\f 12,
\ S_7=-\f{17}8,\ S_9=\f{31}2,\ S_{11}=-\f{691}4,
\\&S_2=S_4=S_6=S_8=S_{10}=0.\endalign$$
 As
$$ (1+\hbox{e}^{-x})\Big(\sum_{n=0}^{\infty}S_n\f{x^n}{n!}\Big)
=\sum_{n=0}^{\infty}\Big(S_n+\sum_{k=0}^n\b nkS_k\Big)\f{x^n}{n!}
=2,$$ we see that
$$\sum_{n=0}^{\infty}S_n\f{x^n}{n!}
=\f{2\hbox{e}^x}{\hbox{e}^x+1}\q(|x|<2\pi).\tag 3.2$$ Since
$\sum_{n=1}^{\infty} S_n\f{x^n}{n!}=\f{\hbox{e}^x-1}{\hbox{e}^x+1}
=-\f{\hbox{e}^{-x}-1}{\hbox{e}^{-x}+1}$, we have $S_n=-(-1)^nS_n$
and so $S_{2m}=0$ for $m\ge 1$. As
$\hbox{e}^{-x/2}\sum_{n=0}^{\infty}S_n\f{x^n}{n!}=\f
2{\hbox{e}^{x/2}+\hbox{e}^{-x/2}}$ is an even function, we have
$$\{S_n\}\in S^+.\tag 3.3$$
Observe that $\sum_{n=0}^{\infty}(\sum_{k=0}^n\b nkE_k)\f{x^n}{n!}
=\hbox{e}^x\cdot\f{2\hbox{e}^x}{\hbox{e}^{2x}+1}$. We also have
$$S_n=\f 1{2^n}\sum_{k=0}^n\b nkE_k
=\f 1{2^n}\sum\Sb k=0\\2\mid k\endSb^n\b nkE_k.\tag 3.4$$
\pro{Corollary 3.4} For $n=1,3,5,\ldots$ we have
$$4S_n+3\sum_{k=1}^{[n/6]}\b n{6k}S_{n-6k}=\cases
2&\t{if $3\nmid n$,}\\-1&\t{if $3\mid n$}.\endcases$$
\endpro
Proof. As $\{S_n\}\in S^+$, from Theorem 3.1 we see that for odd
$n$,
$$3S_n+3\sum\Sb k=1\\3\mid k\endSb^n\b nkS_{n-k}
=3\sum\Sb k=0\\3\mid k\endSb^n\b nkS_{n-k}=\sum_{k=0}^n\b
nkS_k=2-S_n.$$ To see the result, we recall that $S_0=1$ and
$S_{2m}=0$ for $m\ge 1$.

\section*{4. Congruences involving even and odd sequences}
\pro{Theorem 4.1} Let $p$ be an odd prime, and let $\{A_k\}$ be an
odd sequence of rational $p$-integers. Then
$$\sum_{k=1}^{p-1}\f{A_k}{p+k}\e 0\mod{p^2}.$$
\endpro
Proof. Taking $n=p-1$ in Theorem 2.3 we get
$$\sum_{k=0}^{p-1}\b{p-1}k\b{p-1+k}k(-1)^kA_k=0.$$
For $k=1,2,\ldots,p-1$ we see that
$$\align \b{p-1}k\b{p-1+k}k&
=\f{(p-1)(p-2)\cdots(p-k)}{k!}\cdot \f{p(p+1)\cdots(p+k-1)}{k!}
\\&=\f p{p+k}\cdot \f{(p^2-1^2)(p^2-2^2)\cdots(p^2-k^2)}{k!^2}
\\&\e (-1)^k\f p{p+k}\mod {p^3}.\endalign$$
Since $A_0=-A_0$ we have $A_0=0$. Now, from all the above we deduce
the result.
\par\q
\par{\bf Remark 4.1} For given odd prime $p$ and odd sequence
$\{A_n\}$ of rational p-integers, the congruence
$\sum_{k=1}^{p-1}\f{A_k}k\e 0\mod p$ was given by Tauraso[10]
earlier.
 \pro{Corollary 4.1} Let $p$ be an odd prime, and let $\{A_k\}$
be an odd sequence of rational $p$-integers. Then
$$\sum_{k=1}^{p-1}\f{A_k}k\e
p\sum_{k=1}^{p-1}\f{A_k}{k^2}\mod{p^2}.$$
\endpro
Proof. For $k=1,2,\ldots,p-1$ we have $\f 1{k+p}=\f{k-p}{k^2-p^2}\e
\f{k-p}{k^2}=\f 1k-\f p{k^2}\mod{p^2}$. Thus, the result follows
from Theorem 4.1.
 \pro{Corollary 4.2} Let $p$ be an odd prime. Then
$$\sum_{k=1}^{(p-1)/2}\f{(2^{2k}-1)B_{2k}}{p+2k}\e \f{p-1}2\mod{p^2}
\qtq{and}\sum_{k=1}^{(p-1)/2}\f{E_{2k}-1}{(p+2k)2^{2k}}\e 0\mod
{p^2}.$$
\endpro
Proof. Since $\{(-1)^n(2^n-1)B_n\}$ and $\{\f{E_n-1}{2^n}\}$ are odd
sequences, $E_{2m-1}=0$ and $B_{2m+1}=0$ for $m\ge 1$, the result
follows from Theorem 4.1.
 \pro{Corollary 4.3} Let $p$ be an odd
prime, $b,c\in\Bbb Z_p$ and $b(b^2-4c)\not\e 0\mod p$.  Then
$$\sum_{k=0}^{p-1}\f{U_k(b,c)}{(p+k)b^k}\e 0\mod{p^2}.$$
\endpro
Proof. Since $\{\f{U_n(b,c)}{b^n}\}$ is an odd sequence, the result
follows from Theorem 4.1.
\par\q
\par Let $F_n=U_n(1,-1)$ and $L_n=V_n(1,-1)$ be the Fibonacci
sequence and Lucas sequence, respectively. From Section 1 we know
that $\{F_n\}$ is an odd sequence and $\{L_n\}$ is an even sequence.
 \pro{Corollary 4.4} Let $p>5$ be a prime. Then
$$\sum_{k=1}^{p-1}\f{F_k}k\e -\Ls p5\f{5p}4
\Ls{F_{p-\sls p5}}p^2\mod{p^2}.$$
\endpro
Proof. Recently Hao Pan and Zhi-Wei Sun ([4]) proved that
$$\sum_{k=1}^{p-1}\f{F_k}{k^2}\e -\f 15\Ls p5\Ls{L_p-1}p^2\mod p.$$
It is known ([6]) that $F_{p-\sls p5}\e 0\mod p$ and  $L_{p-\sls
p5}\e 2\sls p5\mod{p^2}$. Also, $5F_n=2L_{n+1}-L_n=L_n+2L_{n-1}$.
Thus
$$5F_{p-\sls p5}=2L_p-\Ls p5L_{p-\sls p5}\e 2(L_p-1)\mod{p^2}$$ and
so
$$\sum_{k=1}^{p-1}\f{F_k}{k^2}\e -\f 15\Ls p5\Ls{5F_{p-\sls
p5}}{2p}^2\mod{p^2}.$$ Since $\{F_k\}$ is an odd sequence, applying
Corollary 4.1 we deduce the result.

\pro{Theorem 4.2} Let $p$ be a prime greater than $3$, and let
$\{A_k\}$ be an even sequence. Suppose that
$A_0,A_1,\ldots,A_{p-2},A_p,pA_{p-1}\in\Bbb Z_p$. Then
$$\sum_{k=1}^{p-2}\f{A_k}{p-k}\e \f{2A_p-A_0-pA_{p-1}}p\mod{p^2}.$$
\endpro
Proof. Taking $n=p$ in Theorem 2.3 we get
$$\sum_{k=0}^p\b pk\b{p+k}k(-1)^kA_k=0.$$
For $k=1,2,\ldots,p-1$ we see that
$$\align \b pk\b{p+k}k&
=\f{p(p-1)\cdots(p-k+1)}{k!}\cdot \f{(p+1)\cdots(p+k)}{k!}
\\&=\f p{p-k}\cdot \f{(p^2-1^2)(p^2-2^2)\cdots(p^2-k^2)}{k!^2}
\\&\e (-1)^k\f p{p-k}\mod {p^3}.\endalign$$
Thus,
$$\align &A_0-\b{2p}pA_p+\b
p{p-1}\b{2p-1}{p-1}A_{p-1}+\sum_{k=1}^{p-2}\f p{p-k}A_k
\\&\e\sum_{k=0}^p\b pk\b{p+k}k(-1)^kA_k=0\mod{p^3}.
\endalign$$
Hence
$$\sum_{k=1}^{p-2}\f{A_k}{p-k}
\e\f{2\b{2p-1}{p-1}A_p-p\b{2p-1}{p-1}A_{p-1}-A_0}p\mod{p^2}.$$
 The
famous Wolstenholme's congruence ([13]) states  that
$\b{2p-1}{p-1}\e 1\mod{p^3}.$ Thus the result follows.
 \pro{Corollary 4.5} Let $p$ be a prime
greater than $3$. Then
$$\sum_{k=1}^{(p-3)/2}\f{B_{2k}}{p-2k}\e
\f{p+1}2-\f{pB_{p-1}+1}p\mod{p^2}.$$
\endpro
Proof. It is well known that
$B_0,B_1,\ldots,B_{p-2},B_p,pB_{p-1}\in\Bbb Z_p$. Taking
$A_k=(-1)^kB_k$ in Theorem 4.2 and applying the fact $B_{2k+1}=0$
for $k\ge 1$ we deduce the result.

 \pro{Corollary 4.6} Let $p$ be a prime
greater than $3$, and let $\{A_k\}$ be an even sequence with
$A_0,A_1,\ldots,A_{p-2},A_p,pA_{p-1}\in\Bbb Z_p$. Then
$$\sum_{k=1}^{p-2}\f{A_k}k\e -p\sum_{k=1}^{p-2}\f{A_k}{k^2}+
\f{A_0+pA_{p-1}-2A_p}p\mod{p^2}.$$
\endpro
Proof. For $k=1,2,\ldots,p-2$ we have $\f 1{k-p}=\f{k+p}{k^2-p^2}\e
\f{k+p}{k^2}=\f 1k+\f p{k^2}\mod{p^2}$. Thus, by Theorem 4.2 we
obtain the result.

\pro{Corollary 4.7} Let $p>3$ be a prime, $b,c\in\Bbb Z_p$ and
$b(b^2-4c)\not\e 0\mod p$. Then
$$\sum_{k=1}^{p-1}\f{V_k(b,c)}{(p-k)b^k}\e \f{2(V_p(b,c)-b^p)}{pb^p}\mod{p^2}.$$
\endpro
Proof. Taking $A_k=V_k(b,c)/b^k$ in Theorem 4.2 we deduce the
result.
 \pro{Corollary 4.8} Let $p>5$ be a prime. Then
$$\sum_{k=1}^{p-1}\f{L_k}k\e \f{2(1-L_p)}p\mod{p^2}.$$
\endpro
Proof. Recently Hao Pan and Zhi-Wei Sun ([4]) proved the following
conjecture of Tauraso: $\sum_{k=1}^{p-1}\f{L_k}{k^2}\e 0\mod p.$
Thus taking $A_k=L_k$ in Corollary 4.6 we see that
$$\align\sum_{k=1}^{p-2}\f{L_k}k&\e
-p\Big(\sum_{k=1}^{p-1}\f{L_k}{k^2}-\f{L_{p-1}}{(p-1)^2}\Big)+\f{2+pL_{p-1}-2L_p}p
\\&\e (p+1)L_{p-1}+\f{2(1-L_p)}p\e -\f{L_{p-1}}{p-1}+\f{2(1-L_p)}p
\mod{p^2}.\endalign$$ This yields the result.
 \pro{Theorem 4.3} Let
$p$ be an odd prime and $a_0,a_1,\ldots,a_{\f{p-1}2}\in\Bbb Z_p$. If
$\{a_n\}$ is an even sequence and $p\e 3\mod 4$, or if $\{a_n\}$ is
an odd sequence and $p\e 1\mod 4$, then
$$\sum_{k=0}^{(p-1)/2}\b{2k}k^2\f{a_k}{16^k}\e
\sum_{k=0}^{(p-1)/2}\b{2k}k^2\f{a_{k+2}-a_{k+1}}{16^k}\e 0\mod
{p^2}$$ and
$$\sum_{k=1}^{(p-1)/2}\b{2k}k^2\f{ka_{k-1}}{16^k}\e
\sum_{k=0}^{(p-1)/2}\b{2k}k^2\f{a_{k+1}}{16^k(k+1)}\e
 0\mod {p^2}.$$
\endpro
Proof. Suppose that $\{a_n\}\in S^{\pm}$. By [8, Corollary 3.1],
$\{a_{n+2}-a_{n+1}\}\in S^{\pm}$.  By Lemma 2.2, $\{na_{n-1}\}\in
S^{\mp}$ and $\{\f{a_{n+1}}{n+1}\}\in S^{\mp}$. From Theorem 2.3 we
see that if $\{A_n\}\in S^+$ and $p\e 3\mod 4$, or if $\{A_n\}\in
S^-$ and $p\e 1\mod 4$, then
$$\sum_{k=0}^{\f{p-1}2}\b{\f{p-1}2}k\b{\f{p-1}2+k}k(-1)^kA_k=0.$$
By [9, Lemma 2.2], for $k=0,1,\ldots,\f{p-1}2$ we have
$\b{\f{p-1}2+k}{2k}\e \f 1{(-16)^k}\b{2k}k\mod {p^2}$. Thus,
$$\b{\f{p-1}2}k\b{\f{p-1}2+k}k=\b{2k}k\b{\f{p-1}2+k}{2k}
\e\f 1{(-16)^k}{\b{2k}k^2}\mod {p^2}$$ and so
 $\sum_{k=0}^{\f{p-1}2}\f 1{16^k}{\b{2k}k^2A_k}\e 0\mod{p^2}$ provided
 that $A_0,A_1,\ldots,A_{\f{p-1}2}\in\Bbb Z_p$.
Now combining all the above we deduce the result.

 \pro{Corollary 4.9} Let $p$ be a prime of the
form $4k+3$. Then
$$\sum_{k=0}^{(p-3)/4}\b{4k}{2k}^2\f{B_{2k}}{16^{2k}}\e -\f 18\mod {p^2}.$$
\endpro
Proof. Since $\{(-1)^nB_n\}$ is an even sequence and $B_{2m+1}=0$
for $m\ge 1$, taking $a_n=(-1)^nB_n$ in Theorem 4.3 we deduce the
result.
 \pro{Corollary 4.10} Let $p$ be a prime of the form $4k+1$
and $p=a^2+b^2$ with $a,b\in\Bbb Z$ and $a\e 1\mod 4$. Then
$$\sum_{k=0}^{(p-1)/4}\b{4k}{2k}^2\f{E_{2k}}{32^{2k}}\e 2a-\f p{2a}\mod {p^2}.$$
\endpro
\par Proof. Since  $\big\{\f{E_n-1}{2^n}\big\}$ is an odd sequence,
 taking $a_k=(E_k-1)/2^k$ in Theorem 4.3 and then applying [9, Theorem 2.2]
 we deduce that
$$\sum_{k=0}^{(p-1)/2}\b{2k}k^2\f {E_k}{32^k}\e
 \sum_{k=0}^{(p-1)/2}\b{2k}k^2\f 1{32^k}\e 2a-\f p{2a}\mod {p^2}.$$
 To see the result, we note that $E_k=0$ for odd $k$.

 \pro{Theorem 4.4} Let $p$ be an odd prime and
 $A_0,A_1,\ldots,A_{\f{p-1}2}\in\Bbb Z_p$. If $\{A_n\}\in S^+$ and
 $p\e 3\mod 4$, or if $\{A_n\}\in S^-$ and $p\e 1\mod 4$, then
 $$\sum_{k=0}^{(p-1)/2}\f{\b{2k}k}{2^k}A_k\e 0\mod p.$$
 \endpro
 Proof. Since $\{\f 1{2^n}\}\in S^+$, by Lemma 3.1(i)  we have
 $$\sum_{k=0}^{(p-1)/2}\b{\f{p-1}2}k\Big((-1)^ka_k-(-1)^{\f{p-1}2-k}
 \sum_{s=0}^k\b ks(-1)^sa_s\Big)\f 1{2^{\f{p-1}2-k}}=0.$$
 Note that $\b{\f{p-1}2}k\e \b{-\f 12}k=\f
 1{(-4)^k}\b{2k}k\mod p$. Taking $a_k=A_k$ in the above we deduce
 the result.

 \pro{Theorem 4.5} Let $p$ be an odd prime and
 $A_0,A_1,\ldots,A_p\in\Bbb Z_p$. If $\{A_n\}$ is an odd sequence,
 then
 $$\sum_{k=0}^{(p-1)/2}\f{\b{2k}k}{4^k}A_{p-1-k}
  \e 0\mod p.$$
If $\{A_n\}$ is an even sequence,
 then
 $$\sum_{k=0}^{(p-1)/2}\f{\b{2k}k(k+1)}{4^k}A_{p-2-k}\e
\sum_{k=0}^{(p-1)/2}\f{\b{2k}k}{4^k\cdot k}A_{p-k}\e
 0\mod p.$$
 \endpro
 Proof. Note that $\b{\f{p-1}2}k\e \b{-\f 12}k=\f
 1{(-4)^k}\b{2k}k\mod p$. Taking $n=\f{p-1}2$ in Theorem 2.2 we
 deduce the result.

\end{document}